\documentclass[10pt]{article}

\usepackage{graphicx}
\usepackage{mathptmx}      
\usepackage{graphicx}

\usepackage{latexsym}

\usepackage{amsmath}
\usepackage{amssymb}
\usepackage{amsfonts}

\newcommand{\rf}[1]{(\ref{#1})}

\newcommand{\ba}{\begin{array}}
\newcommand{\ea}{\end{array}}

\newcommand{\be}{\begin{equation}}
\newcommand{\ee}{\end{equation}}

\newcommand{\const}{{\rm const}}

\newcommand{\R}{{\mathbb R}}

\newcommand{\dis}{\displaystyle }

\newcommand{\no}{\noindent}

\newtheorem{prop}{Proposition}[section]
\newtheorem{Th}[prop]{Theorem}

\newtheorem{rem}[prop]{Remark}

\newtheorem{Def}[prop]{Definition}

\newtheorem{ex}[prop]{Example}

\newenvironment{Proof}{\par \vspace{2ex} \par
\noindent \small {\it Proof:}}{\hfill $\Box$ 
\vspace{2ex} \par }

\begin{document}

\title{Energy-preserving numerical schemes of high accuracy for one-dimensional Hamiltonian systems\thanks{This  research work has been partially supported by the grant No.\ N N202 238637 from the Polish Ministry of Science and Higher Education.} 
}

\author{
 {\bf Jan L.\ Cie\'sli\'nski}\thanks{\footnotesize
 e-mail: \tt janek\,@\,alpha.uwb.edu.pl}, \quad 
 {\bf Bogus{\l}aw Ratkiewicz}\thanks{\footnotesize  e-mail: \tt bograt\,@\,poczta.onet.pl} \thanks{\footnotesize permanent address:  I LO, ul.\ \'Sr\'odmie\'scie 31, 16-300 August\'ow, Poland;} 
\\ {\footnotesize Uniwersytet w Bia{\l}ymstoku,
Wydzia{\l} Fizyki, ul.\ Lipowa 41, 15-424
Bia{\l}ystok, Poland}
}

\date{}

\maketitle

\begin{abstract}
We present a class of non-standard numerical schemes which are modifications of the discrete gradient method. They preserve the energy integral exactly (up to the round-off error). The considered class contains locally exact discrete gradient schemes and integrators of arbitrary high order. In numerical experiments we compare our integrators with some other numerical schemes, including the standard discrete gradient method, the leap-frog scheme and a symplectic scheme of 4th order. We study the error accumulation for very long time and the conservation of the energy integral.  
\end{abstract}



\section{Introduction}
\label{intro}

Geometric numerical integration consists in preserving geometric, structural and phy\-sical properties of the considered differential equations. Our aim is to improve the accuracy of geometric integrators, modifying them in an appropriate way, without losing their excellent qualitative properties (including the long-time behaviour, stability and the energy conservation). In this paper we focus on the discrete gradient scheme for one-dimensional Hamiltonian systems. 

Discrete gradient numerical schemes have been introduced many years ago in order to integrate numerically $N$-body systems of classical mechanics with possible applications in molecular dynamics and celestial mechanics \cite{LaG} (see also \cite{Gon,HLW,IA,STW}). Discrete gradient schemes preserve exactly (up to round-off errors) both the total energy and angular momentum. More recently discrete gradient methods have been extended and developed in the context of geometric numerical integration \cite{MQ}. Quispel and his coworkers constructed numerical integrators preserving all integrals of motion of a given system of ordinary differential equations \cite{MQR1,MQR2,QC,QT}. 

In general, geometric numerical integrators are very good in preserving qualitative features of simulated differential equations but it is not easy to enhance their accuracy. Symplectic algorithms can be improved using appropriate splitting methods \cite{Bl,MQ-split,Yo}. Our research is concentrated on improving  the efficiency of the discrete gradient method (which is not symplectic).  

In this paper, we continue our earlier research \cite{CR-long,CR-PRE,CR-PRE2}, extending the theoretical framework on arbitrary one-dimensional Hamiltonian systems. Numerical experiments are carried out in the case of the simple pendulum equation for extremaly long times (we test even 100 millions of periods). We study the accuracy of our new methods, namely, the accumulation of the global error and conservation of the energy integral.

\section{Non-standard discrete gradient schemes} 

In this paper we confine ourselves to one-dimensional Hamiltonian systems
\be  \label{Ham} 
  \dot x = H_p \ , \qquad \dot p = - H_x \ , 
\ee
where $H = H (x, p)$ is a given function, subscripts denote partial differentiation and the dot denotes the total derivative respect to $t$. The Hamiltonian $H (x, p)$ is an integral of motion (the energy integral). 

We consider the following class of non-standard (compare \cite{Mic}) discrete gradients schemes. 
\be  \ba{l} \dis \label{grad-delH} 
 \frac{x_{n+1} - x_n}{\delta_n} =  \frac{  H (x_{n+1}, p_{n+1}) + H (x_n, p_{n+1}) - H (x_{n+1}, p_n) -  H (x_n, p_n) }{2 (p_{n+1}- p_n) } \ , 
\\[4ex] \dis
\frac{p_{n+1} - p_n}{\delta_n} =  \frac{  H (x_n, p_{n+1}) + H (x_n, p_n) - H (x_{n+1}, p_{n+1}) -  H (x_{n+1}, p_n) }{2 (x_{n+1}- x_n) } \ ,   
\ea \ee
where $\delta_n$ is an arbitrary positive function of $h, x_n, p_n, x_{n+1}, p_{n+1}$ etc. (the time step is denoted by $h$). The subscript $n$ indicates that $\delta_n$ may depend on the step $n$.   
In the separable case, i.e., $H = T (p) + V (x)$, the scheme \rf{grad-delH} becomes  
\be  \ba{l} \dis  \label{grad-delTV} 
\frac{x_{n+1} - x_n}{\delta_n} =  \frac{  T  (p_{n+1}) - T (p_n) }{p_{n+1}- p_n } \ ,  
\\[4ex] \dis
\frac{p_{n+1} - p_n}{\delta_n} =  - \frac{  V (x_{n+1}) - V (x_n)}{ x_{n+1}- x_n } \ .
\ea \ee
In numerical experiments we mostly test the case $T(p) = \frac{1}{2}p^2$, where further simplification occurs, see \rf{grad-delV}. 

The system \rf{grad-delH} is a consistent approximation of \rf{Ham}  if we add the condition  
\be  \label{consist} 
 \lim_{h \rightarrow 0} \frac{\delta_n}{h} = 1  \ . 
\ee 
The case $\delta_n = h$ yields the standard discrete gradient method (GR), \cite{Gon,IA,LaG,STW}. 

\begin{Th}  \label{Th-energy} 
The numerical scheme \rf{grad-delH} preserves the  energy integral exactly (up to round-off error), i.e., $H (x_{n+1}, p_{n+1}) = H (x_n, p_n)$. 
\end{Th}

\begin{Proof} The system \rf{grad-delH} implies the equality of both numerators on the right-hand sides of equations \rf{grad-delH}. This, in turn, yields the theorem immediately. 
\end{Proof}

Therefore, any $\delta_n$ satisfying  non-restrictive condition \rf{consist}  yields an energy-pre\-serving numerical scheme. 
The main idea of this paper consists in finding $\delta_n$ such that the resulting numerical scheme is better than the standard gradient method. We consider and test two possibilities. First, the so called locally exact discretizations (section \ref{sec-lex}), then we show that the class \rf{grad-delH} contains integrators of arbitrary high order. The corresponding $h$-series for $\delta_n$ is defined in a recurrent way (section \ref{sec-grad-N}). 

\section{Exact discretization}

 We consider an ordinary differential equation (ODE) with a general solution ${\bf x} (t)$ (satisfying the initial condition ${\bf x} (t_0) = {\bf x}_0$), and a difference equation with the  general solution ${\bf x}_n$. 
The difference equation 
 is the exact discretization of  the ODE  if \ ${\bf x}_n = {\bf x} (t_n)$.

It is well known that any linear ODE with constant coefficients admits the exact discretization in an explicit form \cite{Po}, see also \cite{Ag,CR-ade,Mic}. We summarize these results as follows. 

\begin{Th} \label{Th-exact} 
Any linear equation with constant coefficients, represented in the matrix form by
\be  \label{genlin}
  \frac{d  {\bf x}}{d t} = A {\bf x} + {\bf b} \ ,
\ee
(where ${\bf x} = {\bf x} (t) \in \R^n$, ${\bf b} = \const \in \R^n$ and $A$ is a constant $n\times n$ matrix) admits the exact discretization given by
\be   \label{genex}
   {\bf x}_{n+1} = e^{h A} {\bf x}_n + \left( e^{h A} - I \right) A^{-1} {\bf b}  \ , 
\ee
where \ $h = t_{n+1} - t_n$ \  is the time step and $I$ is the identity matrix. 
\end{Th}

\begin{Proof} The general solution of \rf{genlin} is given by 
\[
{\bf x} (t) = e^{(t-t_0) A} \left( {\bf x} (t_0) + A^{-1} {\bf b} \right) - A^{-1} {\bf b} \ . 
\] 
 Taking into account that that ${\bf x}_n = {\bf x} (t_n)$ and, in particular, ${\bf x}_0 = {\bf x} (t_0)$, we get  
\[
{\bf x}_n  = e^{(t_n-t_0) A} \left( {\bf x} (t_0) + A^{-1} {\bf b} \right) - A^{-1} {\bf b} \ ,
\]
\[
{\bf x}_{n+1}  = e^{(t_n-t_0 + h ) A} \left( {\bf x} (t_0) + A^{-1} {\bf b} \right) - A^{-1} {\bf b} = e^{h  A} \left( {\bf x}_n + A^{-1} {\bf b} \right) - A^{-1} {\bf b}  ,
\]
which ends the proof.   

\end{Proof}

\begin{ex} Exponential growth equation:  
$\dot x = a x$. 
 Exact discretization:
 $x_{n+1} = e^{a h} x_n $ (a geometric series). Equivalent form: 
\be
  \frac{x_{n+1} - x_{n}}{\delta (h) } = a x_n  \ , \qquad  
\delta (h) = \frac{e^{a h} - 1}{a} \ . 
\ee
Note that \  $\dis  \lim_{h\rightarrow 0} \frac{\delta (h) }{h}  =   1 $. 
\end{ex}

\begin{ex}
 Harmonic oscillator: \ 
$\ddot x + \omega^2 x = 0$, \  $p = \dot x$.  
 Exact discretization:  
\be
 x_{n+1} - 2 \cos(\omega h ) x_n + x_{n-1} = 0  , \qquad 
\dis p_n = \frac{ x_{n+1} - \cos(\omega h) x_n }{\sin(\omega h)} \ . 
\ee
\no Equivalent form: 
\be
\frac{x_{n+1} - 2 x_n + x_{n-1}}{ \delta^2 (h) } + \omega^2 x_n + 0 \ , \qquad  \delta (h) =  \frac{2}{\omega} \sin \frac{\omega h}{2} \ . 
\ee
 Note that \ $\delta (h) \approx h$ for $h \approx 0$. 
\end{ex}

Exact discretization seems to be of limited value because, in order to apply it, we need to know the explicit solution of the considered system. However, there exist  non-trivial applications of exact discretizations.  
In the case of the classical Kepler problem we succeeded to use the exact discretization of the harmonic oscillator in two different ways, obtaining numerical integrators preserving all integrals of motion and trajectories \cite{Ci-Kep,Ci-oscyl,Ci-Koz}. Another fruitful direction is associated with the so called locally exact discretizations \cite{Ci-oscyl,Ci-aml,CR-PRE}, see the next section.

\section{Locally exact discrete gradient schemes}
\label{sec-lex}

First, we recall our earlier results concerning the case $H = \frac{1}{2} p^2 + V (x)$, see \cite{CR-long,CR-PRE}. We tested the following class of numerical integrators
\be \ba{l} \label{grad-delV}   \dis
\frac{x_{n+1} - x_n}{\delta_n} = \frac{1}{2} \left( p_{n+1} + p_n \right) \ . \\[3ex] \dis
\frac{p_{n+1} - p_n}{\delta_n} =  - \frac{ V (x_{n+1}) - 
V (x_n) }{x_{n+1} - x_n} \ ,
\ea \ee
where $\delta_n$ is a function defined by
\be \label{deltan} \dis
\delta_n =
\frac{2}{ \omega_n } \tan\frac{ h \omega_n  }{2} \ , \qquad 
\omega_n = \sqrt{   V'' (\bar x)  } \ \ ,  
\ee
and, in general,  $\bar x$ may depend on $n$. For simplicity, we formally assume $V'' (\bar x) > 0$. However, in the case of non-positive $V'' (\bar x) $ one can use the same formula (either, for $V'' (\bar x)  < 0$,  the imaginary unit cancels, or, for  $V'' (\bar x)  = 0$,  we compute the limit $\omega_n \rightarrow 0$ obtaining $\delta_n = h$), for details and final results see \cite{CR-PRE}.

 The simplest choice is $\bar x = x_0$, where $V'(x_0)=0$  (small oscillations around the stable equilibrium). In this case $\delta_n$ does not depend on $n$. The resulting scheme was first presented in \cite{CR-long}, here we propose to name it MOD-GR. In  \cite{CR-PRE} we considered the case 
 $\bar x = x_n$ (which will be called GR-LEX) and its symmetric (time-reversible) modification \ $\bar x = \frac{1}{2} (x_n + x_{n+1})$ \ (GR-SLEX). In both cases  $\bar x$ is changed at every step. 

\begin{Def}
A numerical scheme ${\bf x}_{n+1} = \Psi ({\bf x}_n, h)$ for an autonomous equation $\dot {\bf x} = F ({\bf x})$  is {\it locally exact} if its linearization around any fixed $\bar {\bf x}$ is identical with the exact discretization of the differential equation linearized around $\bar {\bf x}$. 
\end{Def}

We use local exactness as a criterion to select numerical schmemes of high accuracy from a family of non-standard integrators, e.g. from \rf{grad-delH}. Our working algorithm to derive such ``locally exact modifications'' of numerical integrators of the form \rf{grad-delH} assumes that $\delta_n$ depends only on \ $\bar x$, $\bar p$ (or, in more general case,  on $\bar {\bf x}$) and $h$. 
The following theorem extends results of \cite{Ci-oscyl,CR-PRE} on the case of the general time-independent Hamiltonian $H = H (x, p)$. 

\begin{Th}
The discrete gradient scheme \rf{grad-delH} with 
\be  \label{delta-omH} 
  \delta_n = \frac{2}{\omega_n} \tan \frac{h \omega_n}{2} \ , \qquad \omega_n = \sqrt{ H_{xx} H_{pp} - H_{xp}^2 } \ , 
\ee
(where $\omega_n$ is evaluated at $\bar x, \bar p$) is locally exact. 
\end{Th}

\begin{Proof} We have to linearize the continuous system \rf{Ham}, then to find the exact discretization of the obtained linearization.   
Therefore, we put $x = {\bar x} + \xi$, $p = {\bar p} + \eta$ into \rf{Ham} and neglect all terms of order greater than 2. Thus we get
\be  \label{linear-H}
    {\dot \xi} = H_p + H_{px} \xi  + H_{pp} \eta \ , \qquad  {\dot \eta} = - H_x -  H_{xx} \xi - H_{xp} \eta \ . 
\ee
The exact discretization of the system \rf{linear-H} is given by 
\be  \label{dis-Ab} 
  \left( \ba{c} \xi_{n+1} \\ \eta_{n+1} \ea \right)  = e^{h A} \left( \ba{c} \xi_n \\ \eta_n \ea \right)   +  \left( e^{h A} - I \right) A^{-1} {\bf b} \ ,
\ee 
(compare Theorem~\ref{Th-exact}), where
\be  \label{def-Ab} 
  A =  \left( \ba{cc}  H_{xp} & H_{pp} \\ - H_{xx} & - H_{xp} \ea \right) \ , \qquad {\bf b} = \left( \ba{c} H_p \\ - H_x \ea \right) \ .
\ee
We proceed to the linearization of the discrete system \rf{grad-delH}. 
We substitute 
\be
 x_n = {\bar x} + \xi_n \  , \qquad p_n = {\bar p} + \eta_n \ , 
\ee
and assume that $\delta_n$ depends only on ${\bar x}, {\bar p}$ and $h$ (which is equivalent to taking only the first, constant, term of the Taylor expansion of $\delta_n$ with respect to  $\xi_n, \eta_n$).  Then, we linearize the system \rf{grad-delH} around $\bar x, \bar p$ (neglecting terms of at least the second order with respect to $\xi_n$ and $\eta_n$), obtaining
\be  \ba{l} \dis
 \frac{\xi_{n+1} - \xi_n}{\delta_n} = H_p + \frac{1}{2} H_{xp} \left( \xi_n + \xi_{n+1} \right) + \frac{1}{2} H_{pp} \left( \eta_n + \eta_{n+1} \right) \ , \\[3ex]  \dis
 \frac{\eta_{n+1} - \eta_n}{\delta_n} = - H_x - \frac{1}{2} H_{xx} \left( \xi_n + \xi_{n+1} \right) - \frac{1}{2} H_{xp} \left( \eta_n + \eta_{n+1} \right) \ ,
\ea \ee
where partial derivatives $H_x, H_p, H_{xx}, H_{xp}$ and $H_{pp}$ are evaluated at ${\bar x}, {\bar p}$. 
After simple algebraic manipulations we rewrite this linear system in the form
\be  \label{dis-Mw}
  \left( \ba{c} \xi_{n+1} \\ \eta_{n+1} \ea \right)  = M \left( \ba{c} \xi_n \\ \eta_n \ea \right)  + {\bf w} \ ,
\ee 
where
\be  \ba{l} \dis 
   M = \frac{1}{1 + \frac{1}{4} \omega_n^2 \delta_n^2 } \left( \ba{cc}  1 + \delta_n H_{xp} - \frac{1}{4} \omega_n^2 \delta_n^2 & \delta_n H_{pp} \\ - \delta_n H_{xx} & 1 - \delta_n H_{xp} - \frac{1}{4} \omega_n^2 \delta_n^2   \ea \right) \ , 
\\[5ex] \dis 
 {\bf w} = \frac{\delta_n}{1 + \frac{1}{4} \omega_n^2 \delta_n^2 } \left( \ba{cc} 1 + \frac{1}{2} \delta_n  H_{xp} & \frac{1}{2} \delta_n  H_{pp} \\ - \frac{1}{2} \delta_n  H_{xx}  & 1 - \frac{1}{2} \delta_n  H_{xp} \ea \right) \left( \ba{c} H_p \\ - H_x \ea \right) \ , 
\ea \ee
and $\omega_n$ is defined by \rf{delta-omH}. 
Taking into account \rf{def-Ab}, we get
\be  \label{eq-Mw1}
 M = \frac{ 1 - \frac{1}{4} \omega_n^2 \delta_n^2  }{1 + \frac{1}{4} \omega_n^2 \delta_n^2 } + \frac{\delta_n A}{1 + \frac{1}{4} \omega_n^2 \delta_n^2 }  \ , \qquad {\bf w} = \frac{\delta_n {\bf b} }{1 + \frac{1}{4} \omega_n^2 \delta_n^2 }   + \frac{\frac{1}{2} \delta_n^2 A {\bf b} }{1 + \frac{1}{4} \omega_n^2 \delta_n^2 }   
\ee
Systems   \rf{dis-Ab} and \rf{dis-Mw}   coincide  if and only if 
\be  \label{eq-Mw2}
  M = e^{h A} \ , \qquad {\bf w} = \left( e^{h A} - I \right) A^{-1} {\bf b} \ . 
\ee
The proof reduces to showing that the system \rf{eq-Mw2}   is identically satisfied if $\delta_n$ is given by \rf{deltan}. 
We easily verify that 
\be \label{A2}
  A^2 = - \omega_n^2 I \ , \qquad \omega_n^2 = H_{xx} H_{pp} - H_{xp}^2 \ .
\ee
Hence
\be  \label{ehA}
  e^{h A} = \cos h \omega_n + \omega_n^{-1} A \sin h \omega_n \ . 
\ee
The second equation of \rf{eq-Mw2} is satisfied for $\delta_n$ of any form. Indeed, using the first equation of \rf{eq-Mw2} and then the first equation of \rf{A2}, we get 
\be
(M-I) A^{-1} {\bf b} = \frac{ ( \delta_n  - \frac{1}{2} \omega_n^2 \delta_n^2 A^{-1}) {\bf b}  }{1 + \frac{1}{4} \omega_n^2 \delta_n^2 } = \frac{ ( \delta_n  + \frac{1}{2}  \delta_n^2 A ) {\bf b}  }{1 + \frac{1}{4} \omega_n^2 \delta_n^2 } = {\bf w} \ .
\ee
Finally, the first equation of \rf{eq-Mw2} is satisfied if and only if
\be \label{omdel} 
 \frac{ 1 - \frac{1}{4} \omega_n^2 \delta_n^2  }{1 + \frac{1}{4} \omega_n^2 \delta_n^2 } = \cos h\omega_n \ , 
\qquad
\frac{\delta_n }{1 + \frac{1}{4} \omega_n^2 \delta_n^2 } = \frac{\sin h\omega_n}{\omega_n} \ , 
\ee
where we took into account \rf{eq-Mw1} and \rf{ehA}. From the first equation we compute
\be
\frac{1}{1 + \frac{1}{4} \omega_n^2 \delta_n^2 } = \frac{1 + \cos h\omega_n}{2} = \cos^2 \frac{h \omega_n}{2} \ ,
\ee
and substituting it into the second equation of \rf{omdel} we get \rf{delta-omH}. 
\end{Proof}

\begin{rem}
Assuming $\bar x = x_n$, $\bar p = p_n$  we get a numerical scheme called GR-LEX, while the choice $\bar x = \frac{1}{2} \left( x_n + x_{n+1} \right)$, $\bar p = \frac{1}{2} \left( p_n + p_{n+1} \right)$ yields another scheme,  named  GR-SLEX. The system \rf{Ham} is symmetric (time-reversible). The numerical scheme GR-SLEX preserves this property, while GR-LEX does not preserve it. 
\end{rem}

The discrete gradient schemes GR and MOD-GR are of second order. Locally exact discrete gradient schemes have higher order: GR-LEX is of 3rd order and GR-SLEX is of 4th order, see \cite{CR-PRE}. 
In the next section we show how to construct discrete gradient schemes of any order.

\section{Discrete gradient schemes of $N$th order}
\label{sec-grad-N}

We consider the family \rf{grad-delH} of non-standard discrete gradient schemes for the Hamiltonian system \rf{Ham}. The family is parameterized by  a single function $\delta_n$ and this function can be expressed by $x_n, p_n, x_{n+1}, p_{n+1}$ as follows,  
\be
  \delta_n = \frac{2 (x_{n+1} - x_n)(p_{n+1}- p_n) }{  H (x_{n+1}, p_{n+1}) + H (x_n, p_{n+1}) - H (x_{n+1}, p_n) -  H (x_n, p_n) } \ . 
\ee
Replacing here $x_{n+1}, p_{n+1}$ by the exact solution $x (t_{n+1}), p (t_{n+1})$ we formally obtain $\delta_n$ corresponding to the exact integrator. In practice, we can replace $x_{n+1}, p_{n+1}$ by truncated Taylor expansions (and truncate the final result).

Therefore, we take Taylor expansions truncated by neglecting terms of order higher than $N$ (see Appendix A, formulae \rf{TayN}) and compute 
\[
 \frac{2 (x_{n+1}^{[N]}  - x_n)(p_{n+1}^{[N]} - p_n) }{  H (x_{n+1}^{[N]}, p_{n+1}^{[N]}) + H (x_n, p_{n+1}^{[N]}) - H (x_{n+1}^{[N]}, p_n) -  H (x_n, p_n) } = \sum_{k=0}^{N} a_k h^k + O (h^{N+1}) \ , 
\]
where coefficients $a_k$ are functions of $x_n$ and $p_n$. Then, truncating the obtained result, we define
\be
 \delta_n^{[N]} = \sum_{k=1}^N a_k (x_n, p_n) h^k \ .
\ee
The first few coefficients reads
\be  \ba{l}
 a_1 = 1 \ , \quad a_2 = 0 \ , \quad a_3 = H_{xx} H_{pp} - H_{xp}^2 - H_x H_{xpp} - H_p H_{xxp} \ ,  \\[2ex]
a_4  = H_x^2 H_{xppp} - H_p^2 H_{xxxp} + H_p H_{pp} H_{xxx} - H_x H_{xx} H_{ppp} - 3 H_p H_{xp} H_{xxp} \\[2ex]
\qquad  + 3 H_x H_{xp} H_{xpp} \ ,
\ea \ee
where all partial derivatives are evaluated at $x = x_n$, $p = p_n$.  
In the separable case, \ $H = T (p) + V (x)$, \ the formulae simplify
\be  \ba{l}
a_1 = 1 \ , \quad a_2 = 0 \ , \quad a_3 = \frac{1}{12} T_{pp} V_{xx} \ ,
\quad 
 a_4 = \frac{1}{24} \left(  T_p T_{pp} V_{3x} - V_x V_{xx} T_{3p} \right) \ ,
\\[2ex]
 a_5  =  \frac{1}{720} \left(   9 V_x^2 V_{xx} T_{4p} + 9 T_p^2 T_{pp} V_{4x}  - 12 V_x V_{3x} T_{pp}^2 - 12 T_p T_{3p} V_{xx}^2  \right. 
\\[2ex] 
\qquad  +   \left.  6 T_{pp}^2 V_{xx}^2 - 16 T_p T_{3p} V_x V_{3x}       \right) 
\ea \ee
where $V_{k x}$ denotes $k$th derivative of $V$ with respect to $x$, etc.  The case $H = \frac{1}{2} p^2 + V (x)$ is discussed in more detail in \cite{CR-PRE2}, where explicit formulae for $\delta_n^{[N]}$ for $N \leqslant 11$ can be found. 

The gradient scheme \rf{grad-delH} with $\delta_n = \delta_n^{[N]}$ \ is called GR-$N$. Its order is at least $N$, sometimes higher (e.g.,  GR-1 is of 2nd order). Actually  GR-1 and GR-2 are identical with GR.

\section{Numerical experiments}

%
\begin{figure*}
  \includegraphics[width=\textwidth]{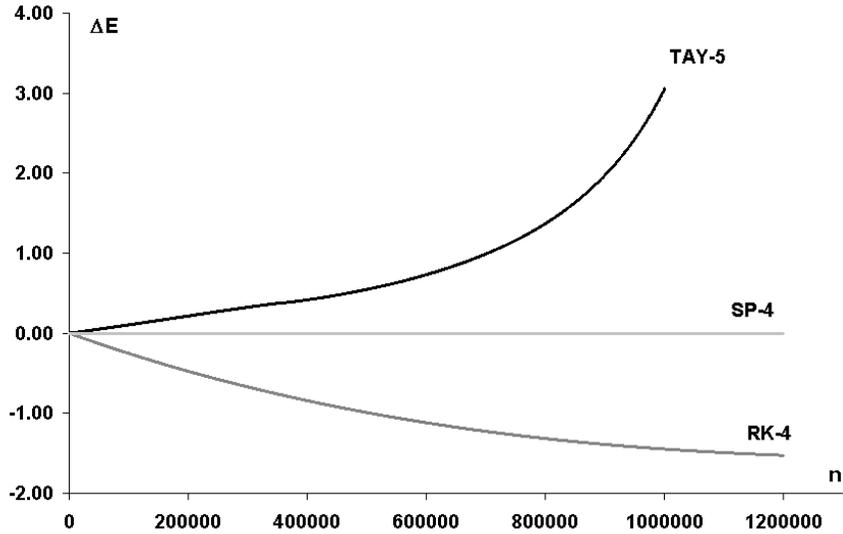}  
\caption{Energy as a function of time ($t = n h$),  $h = 0.25$, $p_0 = 1.8$, $E_{ex} = 0.62$. }
\label{fig-1}       
\end{figure*}
\begin{figure*}[t]
  \includegraphics[width=\textwidth]{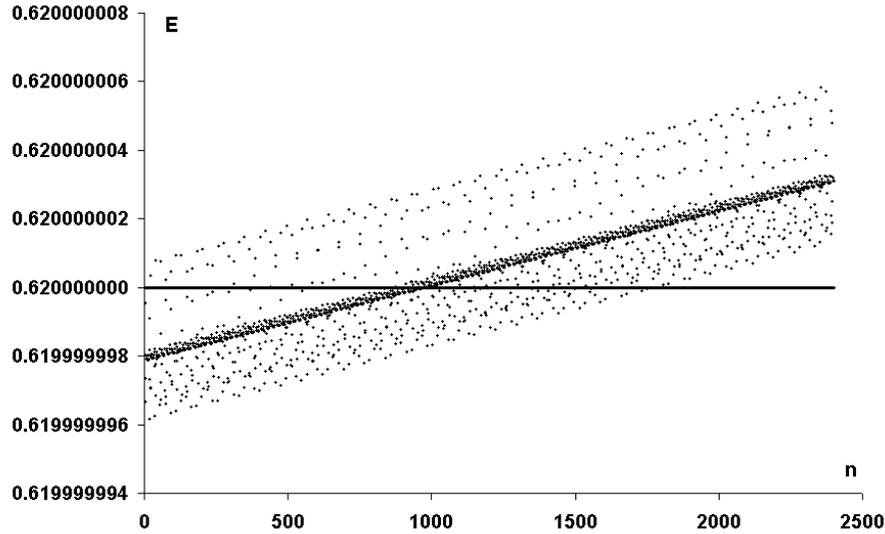}    
\caption{Energy as a function of time ($t = n h$),   $h=0.25$, $p_0 = 1.8$, $E_{ex} = 0.62$. The line $E = 0.62$ corresponds to GR, GR-LEX, GR-3 and GR-7. Other, scattered, points are produced by TAY-10. }
\label{fig-2}       
\end{figure*}

\begin{figure*}
  \includegraphics[width=\textwidth]{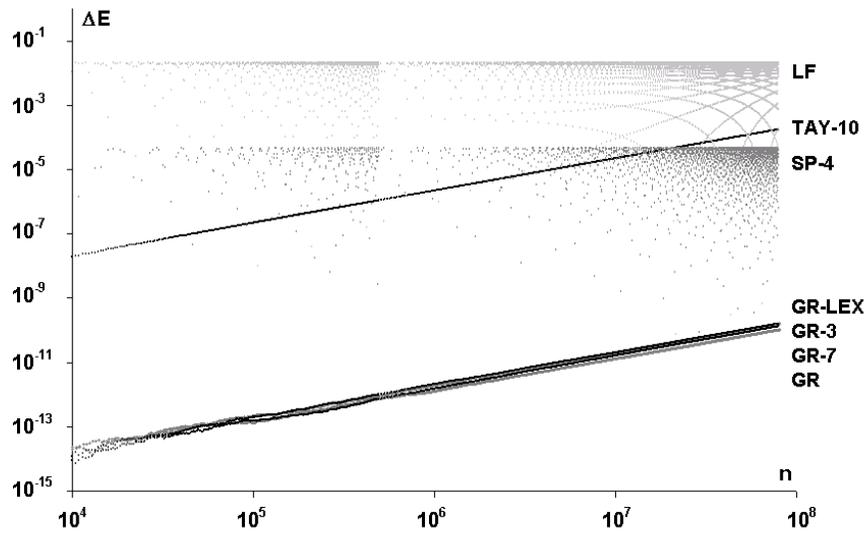}   
\caption{Global error in the energy as a function of time ($t = n h$),   $p_0 = 1.8$, $h = 0.25$. }
\label{fig-3}       
\end{figure*}

\begin{figure*}
  \includegraphics[width=\textwidth]{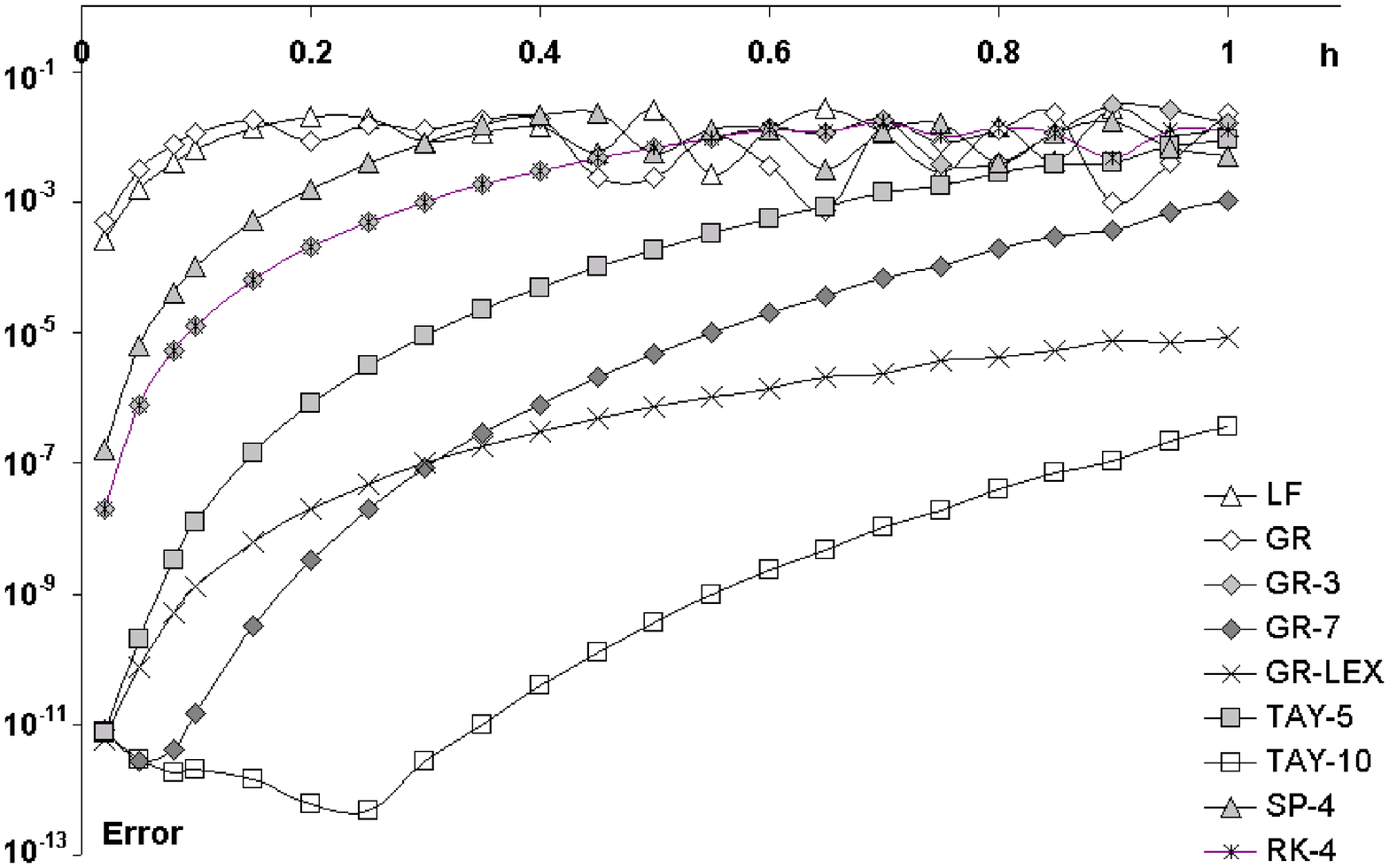}    
\caption{Global error as a function of $h$, evaluated at $t = 120 T_{th}$ for $p_0 = 0.02$. }
\label{fig-4}       
\end{figure*}

\begin{figure*}
  \includegraphics[width=\textwidth]{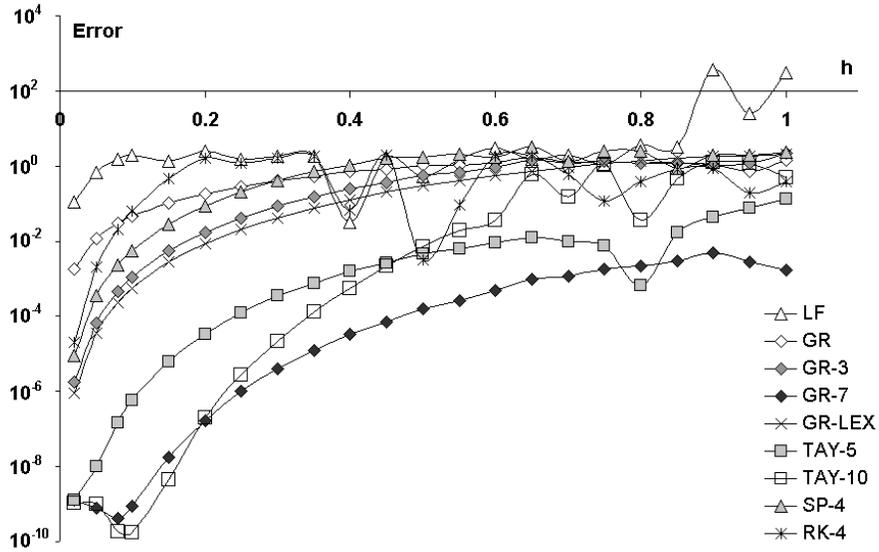}  
\caption{Global error as a function of $h$, evaluated at $t = 120 T_{th}$ for $p_0 = 1.8$.  }
\label{fig-5}       
\end{figure*}

\begin{figure*}
  \includegraphics[width=\textwidth]{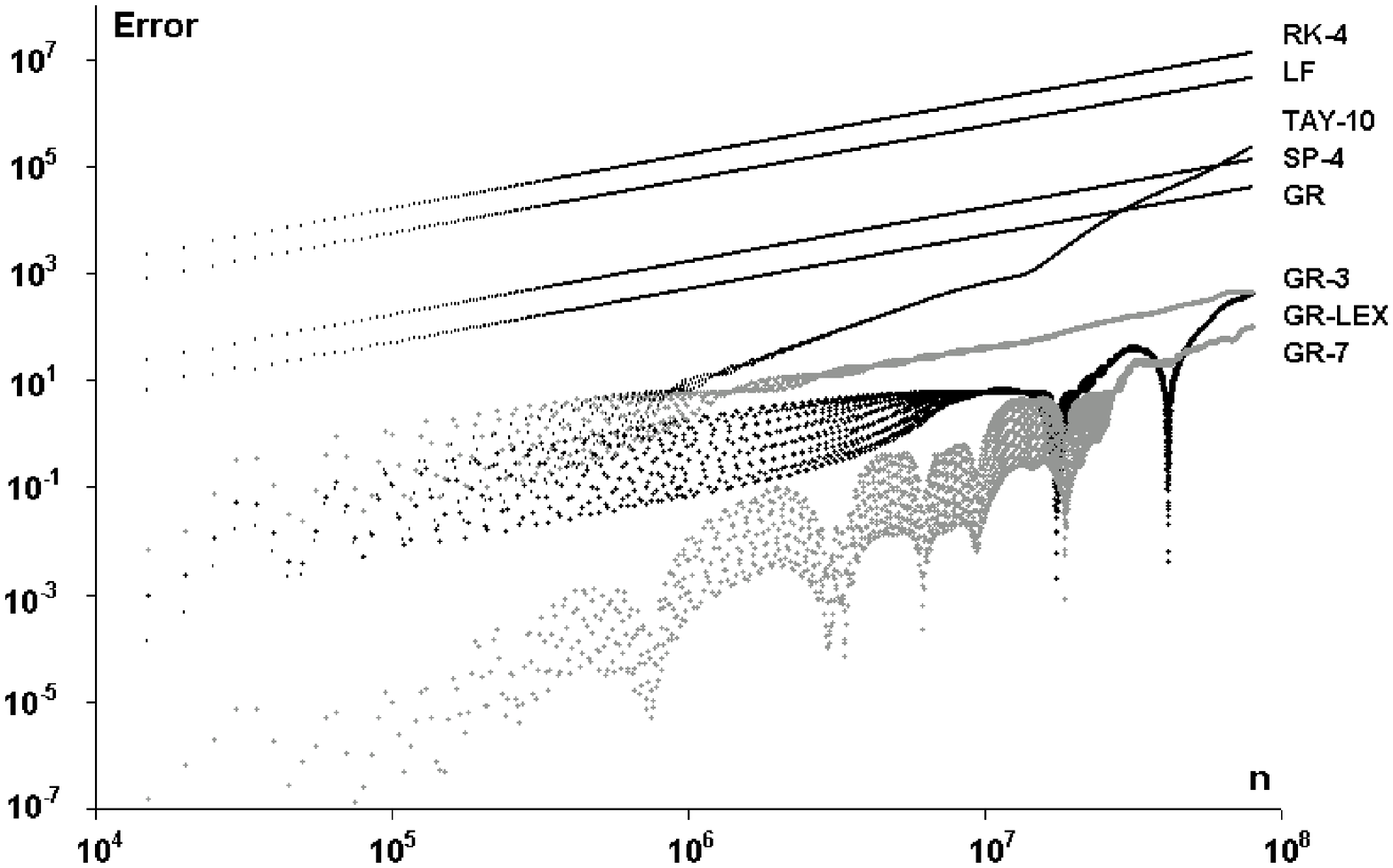}   
\caption{Global error as a function of time ($t = n h$),   $p_0 = 2.001$, $h=0.25$.}
\label{fig-6}       
\end{figure*}

The accuracy of high-order discrete gradient schemes  was tested on the case of the simple pendulum, $H = \frac{1}{2} p^2 - \cos x$ \ (for simplicity always assuming $x_0 = 0$). We compared GR-3, GR-7 and GR-LEX  with the discrete gradient method (GR), the leap-frog scheme (LF), 4th order explicit Runge-Kutta method (RK-4), high-order Taylor methods (TAY-$N$, see appendix A) and a 4th order symplectic scheme (SP-4, see appendix B). Computing global errors we use the exact solution of the simple pendulum equation, expressed in terms of elliptic integrals. 

In previous papers \cite{CR-long,CR-PRE,CR-PRE2} we focused on the stability and accuracy of the period (all motions of the pendulum are periodic). Here, we test the global error, accumulated after $120$ periods (Figures~\ref{fig-4} and \ref{fig-5}) and the accumulation of error after a very long time (up to $n=10^8$ steps), Figure~\ref{fig-6}. We also check the preservation of the energy integral by different numerical schemes, Figures~\ref{fig-1}, \ref{fig-2} and \ref{fig-3}. Details concerning the solution of implicit equations are the same as in \cite{CR-PRE}, e.g., at every step we iterated until the accuracy $10^{-16}$ was  obtained. 

Symplectic methods are known to preserve almost exactly the energy integral \cite{HLW}, some positive results in non-symplectic case are also known \cite{FHP}. 
Figure~\ref{fig-1} shows how accurate is the preservation of the energy by symplectic integrator SP-4 as compared with TAY-5 (permanent, fast growth of the energy) and RK-4 (energy is decreasing approaching the stable equilibrium value). A high order of a given scheme is not sufficient to assure the conservation the energy. From the  beginning TAY-10 produces small, but permanent, drift of the energy, while all discrete gradient schemes yield almost exact value of the energy, see Figure~\ref{fig-2}. 
According to Theorem~\ref{Th-energy} all gradient schemes preserve the energy integral exactly (up to round-off errors). Only after very long time one can notice that also discrete gradient schemes have a slight drift of the energy. A curious phenomenon can be observed at Figure~\ref{fig-3}. Gradient schemes and TAY-10 show a linear growth of the energy error (but the energy error of TAY-10 is always greater by 6 orders of magnitude!), while the energy errors of symplectic schemes, LF and  SP-4,  vary in a large range but do not show any systematic time dependence. 
However, in the considered time interval,  discrete gradient methods preserve the energy  more accurately by several orders of magnitude than symplectic integrators like LF or SP-4.

Figures~\ref{fig-4}, \ref{fig-5} \ show $h$-dependence of the global error calculated at $t = 120 T_{ex}$, where $T_{ex}$ is the period of the exact solution (e.g., $T_{ex} = 6,283342396$
\ for \ $p_0 = 0.02$ 
and $T_{ex} = 9.12219655$ \ for \ $p_0 = 1.8$). We observe that, usually, higher order  integrators are more accurate. An important exception is GR-LEX, of 3rd order, which for $p_0 = 0.02$ and $h > 0.3$ is better than GR-7 (the behaviour of GR-SLEX is almost the same as GR-LEX). However, TAY-10 is clearly the best in this case. 
Only  after very long time evolution TAY-10 becomes less accurate than GR and SP-4, although initially it was comparable with gradient schemes of high order, see Figure~\ref{fig-6}. 
Several schemes at Figure~\ref{fig-6} show linear growth (at least for large $t$). It has been shown, see \cite{CH,HLW}, that symplectic integrators (under some mild conditions) have linear error growth. Results of our experiments suggest that after sufficiently large time some other schemes (e.g., RK-4, GR, GR-3) also accumulate error linearly.  Finally, we point out that until $n=10^7$ the global error of GR-LEX is smaller than the period (and the error of GR-7 is even smaller, by one order of magintude).  The global error of GR is smaller than that of SP-4, not saying about GR-3 or GR-7, see Figure~\ref{fig-6}.

\section{Conclusions}

Modifications presented in this paper \ essentially improve the discrete gradient method (in the one-dimensional case) keeping all its advantages.  Modified gradient schemes GR-LEX, GR-SLEX, GR-$N$ have important advantages: 
\begin{itemize}
\item  conservation of the energy integral (up to round-off errors), 
\item high stability, exact trajectories in the phase space, 
\item high accuracy (third, fourth and $N$th order, respectively), 
\item very good long-time behaviour of numerical solutions.
\end{itemize}
We point out, however, that numerical schemes \rf{grad-delH}, like all discrete gradient methods, are neither symplectic nor volume-preserving.  Most of them, including GR-LEX and GR-$N$ ($N > 2$) are not symmetric (time-reversible). GR and GR-SLEX are symmetric. 

In the near future we plan to generalize the approach presented in this paper 
on some multidimensional cases \cite{Ci-aml} (the crucial point is that $\delta_n$ is a matrix)  and to extend the range of its applications on some other numerical integrators (including the implicit midpoint rule and numerical schemes which preserve integrals of motion \cite{MQR1,QC}). One can also use a variable time step, if needed \cite{Ci-aml}.

\pagebreak

\section*{Appendix A. Explicit Taylor schemes of $N$th order}

$N$th Taylor method for the system  \rf{Ham} is defined by
\be \label{TayN}
   x_{n+1}^{[N]} = \sum_{k=0}^N \frac{b_k h^k}{k!} \ , \qquad p_{n+1}^{[N]} = \sum_{k=0}^N \frac{c_k h^k}{k!} \ , 
\ee
where the coefficients $b_k, c_k$ are computed from Taylor's expansion of the exact solution (see, for instance, \cite{Is-book}, p.18). We assume \ $x (t) = x_n$, $p (t) = p_n$, $t = t_n$ \ and expand $x (t + h)$ and $p (t + h)$ in Taylor series:
\be
   x (t+h) = \sum_{k=0}^\infty \frac{b_k h^k}{k!} \ , \qquad b_k = \left. \frac{d^k x (t)}{d t^k} \right|_{t = t_n} \ ,
\ee
\be
  p (t+h) = \sum_{k=0}^\infty \frac{c_k h^k}{k!}  \ , \qquad  c_k = \left. \frac{d^k p (t)}{d t^k} \right|_{t=t_n} \ ,
\ee
where all derivatives are replaced by functions of $x_n, p_n$ \ using \rf{Ham} and its differential consequences, e.g., 
\be
\ddot x = \frac{d}{d t} \frac{\partial H}{\partial p} = H_{px} {\dot x} + H_{pp} {\dot p} = H_{px} H_p - H_{pp} H_x \ ,
\ee
where $H_x, H_p, H_{xx}$ etc. are evaluated at $x_n, p_n$. Thus  we get 
\be \ba{l}  \label{bck}
b_1 = H_p \ , \qquad b_2 = H_p H_{xp} - H_x  H_{pp} \ , 
\\[2ex]
b_3 = H_x^2 H_{3p} + H_{xxp} H_p^2 - 2 H_x H_p H_{ppx}   + H_p H_{xp}^2 - H_p H_{pp} H_{xx} \ , 
\\[2ex]
c_1 = - H_x \ , \qquad c_2 = H_x H_{xp} - H_p H_{xx} \ , 
\\[2ex]
c_3 = - H_p^2 H_{3x} - H_{xpp} H_x^2 + 2 H_x H_p H_{xxp}  - H_x H_{xp}^2 + H_x H_{pp} H_{xx} \ ,
\ea \ee
and subsequent coefficients can be easily computed using the total derivative: 
\be
  b_{k+1} = \frac{d b_k}{d t} = H_p \frac{\partial b_k}{\partial x} - H_x \frac{\partial b_k}{\partial p} \ ,  \qquad 
c_{k+1} = H_p \frac{\partial c_k}{\partial x} - H_x \frac{\partial c_k}{\partial p} \ . 
\ee
The  formulae \rf{bck} simplify in the case \ $H = T (p) + V (x)$ \cite{CR-PRE2}. 

\pagebreak

\section*{Appendix B. Explicit symplectic schemes of $2 M$th order}

Symplectic explicit integrators of arbitrary even order $N = 2 M$ can be derived by composition methods, see, e.g., \cite{HLW}. In this section we present results of the pioneering paper \cite{Yo}, confining ourselves to the case $H = \frac{1}{2} p^2 + V (x)$. 

The numerical scheme SP-2$M$ is defined by the following procedure. 
Having $x_n, p_n$ we compute the next step, $x_{n+1}, p_{n+1}$, as follows. We denote $x_{n} = x^{[0]}, p_n = p^{[0]}$ and perform $K +1$ iterations (where $K = 3^{M-1}$) 
\be
x^{[i]} = x^{[i-1]} + h c_M^{[i]} p^{[i-1]} \ , \quad 
p^{[i]} = p^{[i-1]} - h d_M^{[i]} V' (x^{[i]}) \ ,
\ee
where $c_M^{[i]}, d_M^{[i]}$ ($i=1,2,\ldots,K+1$) have to  be carefully computed (see below) in order to secure the required order.
Then we identify $x^{[K+1]} = x_{n+1}, p^{[K+1]} = p_{n+1}$. 

The coefficients $c_M^{[i]}, d_M^{[i]}$ are computed recursively. First, all coefficients for $M=1$ are given by: 
\be
  c_1^{[1]} = \frac{1}{2} \ , \quad d_1^{[1]} = 1 \ , \quad c_1^{[2]} = \frac{1}{2} \ , \quad d_1^{[2]} = 0 \ . 
\ee
Then, we express coefficients  $c_{m+1}^{[i]}, d_{m+1}^{[i]}$ by coefficients $c_m^{[i]}, d_m^{[i]}$:
\be \ba{l}
d_{m+1}^{[i]} = d_{m+1}^{[2k+i]} = y_m d_m^{[i]} \ , \quad (i=1,\ldots,k) \ , 
\\[2ex]
d_{m+1}^{[k+i]} = (1-2 y_m) d_m^{[i]} \ ,  \quad (i=1,\ldots,k) \ , 
\\[2ex]
d_{m+1}^{[3k +1]} = 0 \ , 
\\[2ex]
c_{m+1}^{[i]} = c_{m+1}^{[2k+i+1]} = y_m c_m^{[i+1]} \ , \quad (i=1,2,\ldots,k) \ , 
\\[2ex]
c_{m+1}^{[k+i+1]} = (1-2 y_m) c_m^{[i+1]} \ ,  \quad (i=1,2,\ldots,k-1) \ , 
\\[2ex]
c_{m+1}^{[k +1]} = (1-y_m) \left( c_m^{[1]} + c_m^{[k+1]} \right) \ , 
\\[2ex]
c_{m+1}^{[2k +1]} = (1-y_m) \left( c_m^{[1]} + c_m^{[k+1]} \right) \ , 
\ea \ee
where $k = 3^{m-1}$ and   
\be \label{ym} 
  y_m = \frac{1}{2 - 2^{1/(2 m + 1)}} \ .
\ee
In particular, the symplectic integrator SP-4 has the following coefficients
\be \ba{l} \dis
 c_2^{[1]} = c_2^{[4]} = \frac{1}{2 (2 - 2^{1/3})} \ ,  \quad  c_2^{[2]} = c_2^{[3]} = \frac{1 - 2^{1/3}}{2 (2 - 2^{1/3})} \ ,  \\[3ex] \dis
d_2^{[1]} = d_2^{[3]} = \frac{1}{2 - 2^{1/3}} \ ,  \quad d_2^{[2]} = - \frac{2^{1/3}}{2-2^{1/3}} \ , \quad 
d_2^{[4]} =  0 \ .
\ea \ee
The scheme SP-4 was independently presented in \cite{FR} and \cite{Yo}.

\pagebreak



\end{document}